\documentclass[leqno]{article}
\usepackage{verbatim, amsmath, amscd, amssymb}
\usepackage[all]{xy}
\newcommand{\bm}{{\bf m}}
\newcommand{\bU}{{\bf U}}

\newcommand{\cA}{{\mathcal A}}
\newcommand{\cB}{{\mathcal B}}
\newcommand{\cC}{{\mathcal C}}

\newcommand{\cK}{{\mathcal K}}

\newcommand{\cV}{{\mathcal V}}

\newcommand{\fru}{\mathfrak u}

\newcommand{\rA}{\mathrm{A}}

\newcommand{\bbC}{\mathbb C}

\newcommand{\bbF}{\mathbb F}

\newcommand{\bbN}{\mathbb N}
\newcommand{\bbQ}{\mathbb Q}

\newcommand{\bbZ}{\mathbb Z}

\newcommand{\ch}{\mathrm{ch}}

\newcommand{\Dist}{\mathrm{Dist}}

\newcommand{\res}{\mathrm{res}}

\newcommand{\SL}{\mathrm{SL}}

\newcommand{\Mod}{\mathbf{Mod}}

\newcommand{\lbr}{\begin{bmatrix}}
\newcommand{\rbr}{\end{bmatrix}}
\newcommand{\for}{\bigcirc\kern-2.6ex \because}
\newcommand{\forb}{\bigcirc\kern-2.8ex \because}
\newcommand{\forbb}{\bigcirc\kern-3.0ex \because}
\newcommand{\forbbb}{\bigcirc\kern-3.1ex \because}

\newcommand\pf{\noindent {\bf Proof:  }}
\newtheorem{thm}{Theorem.}

\newtheorem{prop}{Proposition.}

\newtheorem{cor}{Corollary.}

\begin{document}
\large
\title{
\bf
On certain maximal
cyclic modules
for the quantized special linear algebra
at a root of unity
\thanks
{supported in part
by JSPS Grants in Aid
for Scientific Research}}

\author{
K\textsc{aneda} Masaharu
\\
558-8585
Osaka, Japan
\\
Osaka City University
\\
Department of
Mathematics
\\
kaneda@sci.osaka-cu.ac.jp
\and
N\textsc{akashima} Toshiki
\\
102-8554 Tokyo, Japan
\\
Sophia University
\\
Department of Mathematics
\\
toshiki@mm.sophia.ac.jp}
\maketitle

\begin{abstract}
By properly specializing the parameters
irreducible modules of maximal dimension for the De Concini-Kac
version of the Drinfeld-Jimbo
quantum algebra in type $A$
may be
transformed into modules over Lusztig's infinitesimal quantum algeba.
Thus obtained modules have a simple socle and a simple head,
and share the same dimension as
the infinitesimal Verma modules.
Despite these common features
we find that they are never isomorphic to infinitesimal Verma modules
unless they are irreducible.
The same carry over to the modular setup for the special linear groups in
positive characteristic.
\end{abstract}

The finite dimensional irreducible representations of
the De Concini-Kac version
of the Drinfeld-Jimbo quantized
enveloping algebra
at a
complex $\ell$-th root of 1
have their dimensions bounded above,
generically attaining a maximal
dimension
\cite{DK}.
In type $\rA$
Date, Jimbo, Miki and Miwa \cite{DJMM}
have given a concrete
realization
of most of those of maximal dimension.
By properly specializing their parametres the second named author of the
present paper
found in \cite{N}
that
they afford modules $\cV$, rarely irreducible,
for Lusztig's
\lq\lq
infinitesimal"
quantum algebra
$\fru$
\cite{L1}
and that each
$\cV$ has a unique, up to scalar, invariant vector
$u_{\vec{0}}$
relative to a Borel subalgebra
$\fru^\sharp$
of
$\fru$,
and hence that
$\cV$ has a simple socle generated by
$u_{\vec{0}}$.
The dimension being right, it is tempting to compare $\cV$
with
Humphreys'
\lq\lq
infinitesimal" Verma modules
\cite{H}
quantized by
Andersen, Polo and Wen
\cite{APW},
which are the standard objects
of study in the representation theory of
$\fru$.

It is easy to see that
$\cV$ is isomorphic to
an infinitesimal Verma module as $\fru^\sharp$-module,
which in turn shows that
$\cV$ has the same simple head as the infinitesimal Verma module.
The explicit description of the actions of the
standard generators of
$\fru$ on
$\cV$ allows us, however, to find that $\cV$ has also a unique,
up to scalar, invariant vector
$u_-$
with respect to
the opposite infinitesimal Borel subalgebra.
It follows that
$\cV$
does not lift to an integrable
$\fru
U^0_\bbC$-module,
and
hence
$\cV$ can not be
isomorphic to any infinitesimal Verma module as $\fru$-module
unless
$\cV$ is simple,
where
$U^0_\bbC$
is the Cartan part of
Lusztig's quantum algebra
$U_\bbC$
at the
$\ell$-th root of 1.

By construction $\cV$ may be defined over
$\cB=\bbZ[v, v^{-1}]/(\phi_\ell)$
with
$\phi_\ell$
the
$\ell$-th cyclotomic polynomial in indeterminate
$v$.
If $\ell$ is an odd prime
$p$,
$\cB/p\cB$
is a finite field
$\bbF_p$ of $p$-elements.
Let
$G$
be the special linear group scheme
over
$\bbF_p$
with opposite Borel subgroups
$B$ and $B^+$,
and let
$G_1$, $B_1$, $B^+_1$
be the Frobenius kernel of
$G$, $B$, and $B^+_1$,
respectively.
If $\cV_\cB$ is the
$\cB$-form of $\cV$,
then
$\cV_p=\cV_\cB\otimes_\cB\bbF_p$
is naturally a
$G_1$-module.
We find that
$u_{\vec{0}}\otimes1$
(resp. $u_-\otimes1$)
remains a unique, up to $\bbF_p^\times$,
$B^+_1$-
(resp. $B_1$-)invariant vector
in
$\cV_p$.
Hence $\cV_p$
is isomorphic to an infinitesimal Verma module as
$B^+_1$-module, but not as $G_1$-module
unless
$\cV_p$ is simple.

If the simple $\fru$-module
generated by
$u_{\vec{0}}$
in
$\cV$ and the simple $G_1$-module generated by $u_{\vec{0}}\otimes1$ in $\
cV_p$
have the same dimension,
Lusztig's conjecture for the irreducible characters of
$G$-modules will follow from
the celebrated theorems of
Kazhdan and Lusztig \cite{KL}
and
Kashiwara and Tanisaki \cite{KT}
and Casian
\cite{C}.

If $\cC$ is a category,
$\cC(A, B)$
will denote the set of morphisms of
$\cC$ from object
$X$ of
$\cC$ to
object $Y$
of $\cC$.
If $A$ is a ring,
$A\Mod$
will denote the category of
left
$A$-modules.

We are grateful to H.H. Andersen
and J.C. Jantzen
for a helpful comment.


\ \\
{\bf
$1^\circ$
Infinitesimal Verma modules}

In this section we recollect some facts about infinitesimal Verma modules
over an arbitrary quantum algebra of finite type.

\ \\
(1.1)
Let
$\bbQ(v)$ be the fractional field of the Laurent polynomial ring
$\cA=\bbZ[v, v^{-1}]$
in indeterminate
$v$,
$A=[\!(A_{ij})\!]$
an indecomposable Cartan matrix of finite type and
let
$\bf U$ be the associated
Drinfeld-Jimbo
quantum algebra over
$\bbQ(v)$
with generators
$E_i$,
$F_i$,
and
$K_i^{\pm1}$,
$i\in [1, n]$.
Let
$U$ be Lusztig's
$\cA$-subalgebra of
$\bU$ generated by
$E_i^{(r)}=
\frac{E_i^r}{[r]_i!}$,
$F_i^{(r)}=
\frac{F_i^r}{[r]_i!}$,
$K_i^{\pm1}$,
$i\in[1, n]$,
$r\in\bbN$,
where
$[r]_i!=
\prod_{s=1}^r[s]_i$
with
$[s]_i=
\frac{v_i^s-v_i^{-s}}{v_i-v_i^{-1}}$,
$v_i=v^{d_i}$,
$d_i\in\{1, 2, 3\}$
minimal such that
the matrix
$[\!(d_iA_{ij})\!]$
is symmetric.

Let
$R$
(resp.
$\Lambda$)
be the root system
(resp. the weight lattice)
associated to
$A$
and
$R^+$
a positive subsystem of
$R$
with the simple roots
$\alpha_i$,
$i\in[1, n]$.
We equip
$\Lambda$ with a partial order
defined by
$R^+$
as usual.
Let
$\Lambda^\vee$ be the colattice of $\Lambda$
and denote by
$\langle\ , \ \rangle : \Lambda\times\Lambda^\vee
\to
\bbZ$ the perfect pairing.
If
$\alpha\in R$,
let $\alpha^\vee$ be its coroot.
We set
$\lambda_i=\langle\lambda,
\alpha_i^\vee\rangle$
for
$\lambda\in \Lambda$.
Let
$U^0$
be the $\cA$-subalgebra of
$U$
generated by
$K_i^{\pm1}$
and
$\lbr
K_i; c
\\
r\rbr=
\prod_{s=1}^r
\frac{K_iv_i^{c-s+1}-K_i^{-1}v_i^{-c+s-1}}
{v_i^s-v_i^{-s}}$,
$i\in[1, n]$,
$c\in\bbZ$,
$r\in\bbN$.
Each
$\lambda\in\Lambda$
defines an
$\cA$-algebra homomorphism
$\chi_\lambda : U^0\to\cA$
such that
\[
K_i\mapsto
v_i^{\lambda_i},
\quad
\lbr
K_i; c
\\
r\rbr
\mapsto
\lbr
\lambda_i+c
\\
r
\rbr_i
=
\frac{\prod_{s=0}^{r-1}
[\lambda_i+c-s]_i}{[r]_i!}
\quad\forall
i\in[1, n],
c\in\bbZ,
r\in\bbN.
\]

\setcounter{equation}{0}
\ \\
(1.2)
Let
$\ell$
be a positive integer greater than 2 prime to all entries
$A_{ij}$
of the Cartan matrix $A$,
$\cK=\bbQ[v]/(\phi_\ell)$,
and set
$U_\cK=U\otimes_\cA\cK$.
Let
$\fru$
(resp. $\fru^+$;
$\fru^-$;
$\fru^0$)
be
the $\cK$-subalgebra of
$U_\cK$
generated by
$E_i\otimes1$,
$F_i\otimes1$,
$K_i\otimes1$
(resp.
$E_i\otimes1$;
$F_i\otimes1$;
$K_i\otimes1$),
$i\in[1, n]$.
Let also
$\fru^\sharp=\fru^+\fru^0$
and
$\fru^\flat=\fru^0\fru^-$.
We will abbreviate
$x\otimes1$ of $U_\cK$
as
$x$,
and
$\chi_\lambda\otimes_\cA\cK$
as
$\chi_\lambda$.

Let
$\tilde\fru$
be the $\cK$-subalgebra of
$U_\cK$
generated by
$\fru$ and
$U^0_\cK=U^0\otimes_\cA\cK$,
and let
$\tilde\fru^\sharp=\fru^\sharp
U^0_\cK$,
$\tilde\fru^\flat=
U^0_\cK\fru^\flat$.
Each
$\lambda\in\Lambda$
defines a
1-dimensional
$\tilde\fru^\flat$-module
by
$\chi_\lambda$
annihilating
all
$F_i$,
which we will still
denote by
$\lambda$.
Let
$\tilde\nabla(\lambda)=
\tilde\fru^\flat\Mod(\tilde\fru,
\lambda)$.
We make
$\tilde\nabla(\lambda)$
into
a
$\tilde\fru$-module
by setting
$
xf=f(?x)$
for each
$
x\in\tilde\fru$ and
$f\in\tilde\nabla(\lambda)$.

Let
$\Lambda^\res=\{
\nu\in\Lambda
\mid
\nu_i\in[0, \ell-1]
\ \forall i
\}$.
If we write
$\lambda=
\lambda^0+\ell\lambda^1$
with
$\lambda^0\in\Lambda^\res$
and
$\lambda^1\in\Lambda$,
one has from
\cite[1.8]{APW}
an isomorphism of
$\tilde\fru$-modules
\begin{equation}
\tilde\nabla(\lambda)\simeq
\tilde\nabla(\lambda^0)\otimes_\cK\ell\lambda^1,
\end{equation}
where
$\ell\lambda^1$
ia a
1-dimensional
$\tilde\fru$-module
defined by
$\chi_{\ell\lambda^1}$
annihilating all
$E_i$,
$F_i$
and
$K_i-1$.

On the other hand,
the natural gradation on
$\fru^+$
assigning
each
$E_i$ grade
$\alpha_i$
equips
$\fru^+$
with a structure of
$\tilde\fru^\sharp$-module
such that $\fru^+$
act by the left multiplication
and
$U^0_\cK$
by
$\chi_{-\lambda+\nu}$
on the $\nu$-th homogeneous part
of
$\fru^+$,
$\nu\in\sum_i\bbN\alpha_i$.
Recall antiautomorphism
$\Psi$ on $U$
such that
$E_i\mapsto
E_i$,
$F_i\mapsto
F_i$
and
$K_i\mapsto
K_i^{-1}$,
$i\in[1, n]$.
If $M$
is a
$\tilde\fru^\sharp$-module of finite type,
we will denote by
$M^\Psi$
the
$\cK$-linear dual of
$M$
made into
$\tilde\fru^\sharp$-module
by setting
$xf=f(\Psi(x)?)$
for each
$x\in\tilde\fru^\sharp$,
$f\in\Mod_\cK(M, \cK)$.
Then we have an isomorphism of
$\tilde\fru^\sharp$-modules
\begin{equation}
\tilde\nabla(\lambda)\simeq
(\fru^+)^\Psi
\quad\text{via}\quad
f\mapsto
f\circ(\Psi\otimes_\cA\cK).
\end{equation}

One can likewise define a
$\fru$-module
$\nabla(\lambda)=
\fru^\flat\Mod(\fru,
\lambda)$.
By restricting the
$\tilde\fru$-action to $\fru$,
$\tilde\nabla(\lambda)$
yields
$\nabla(\lambda)$.
Then
the isomorphism
(2) restricts to an isomorphism of
$\fru^\sharp$-modules
\begin{equation}
\nabla(\lambda)\simeq
(\fru^+)^\Psi.
\end{equation}

\setcounter{equation}{0}
\noindent
(1.3)
Recall from Xi \cite[2.5]{X}
that
\begin{equation}
\text{
$\fru^+$
has a simple socle
$\cK\prod_{\alpha\in R^+}
E_\alpha^{\ell-1}$
as
$\fru^+$-module},
\end{equation}
where
$E_\alpha$ is a root vector of
$\fru^+$
associated to
$\alpha\in R^+$
\cite{L2}
and the product is taken in a
certain
specific order.
It follows that
$\fru^+$
is indecomposable as
$\fru^+$-module,
and
hence that
\begin{equation}
\text{
$\fru^+$
is a projective cover of
trivial module
$\cK$
as
$\fru^+$-module}.
\end{equation}

By an integrable
$\tilde\fru$-
(resp. $\fru$-)
module
$M$
we will mean a
$\tilde\fru$-
(resp. $\fru$-)
module
$M$
such that
\[
M=\coprod_{\mu\in \Lambda}
M_\mu,
\quad
\text{with}
\quad
M_\mu=
\{
m\in M\mid
tm=\chi_\mu(t)m
\ \forall t\in U^0_\cK
\text{ (resp. $\fru^0$)}\}.
\]
Each
$\tilde\nabla(\lambda)$
(resp.
$\nabla(\lambda)$)
is an integrable
$\tilde\fru$-
(resp.
$\fru$-)
module.
Define integrable $\tilde\fru^\sharp$-
and
$\fru^\sharp$-modules
likewise.
One obtains from (2) and
(1.2.2, 3)

\begin{prop}
For each
$\lambda\in\Lambda$
the
$\fru^\sharp$-
(resp. $\tilde\fru^\sharp$-)
module
$\nabla(\lambda)$
(resp.
$\tilde\nabla(\lambda)$)
is an injective hull of
$\lambda$
in the category of integrable
$\fru^\sharp$-
(resp.
$\tilde\fru^\sharp$-)
modules.
\end{prop}

\setcounter{equation}{0}
\noindent
(1.4)
Let
$\zeta$
be the image of
$v$
in
$\cK$,
and let
$U_\zeta$
be the De Concini-Kac
algebra
\cite{DK}
over
$\cK$
associated to the Cartan matrix
$A$
with the generators
$E_i$,
$F_i$,
$K_i^{\pm1}$,
$i\in[1, n]$,
and
the same relations
for
$\bf U$
with
$v$ replaced by
$\zeta$.
For each
$\alpha\in R^+$
let
$E_\alpha$
(resp.
$F_\alpha$)
be the root vector of
$U_\zeta$
associated with
$\alpha$
(resp.
$-\alpha$),
and let
$\overline
U_\zeta=
U_\zeta/
(E_\alpha^\ell,
F_\alpha^\ell
\mid
\alpha\in
R^+$).
If
$\overline
U_\zeta^\flat$
is the
$\cK$-subalgebra
of
$\overline
U_\zeta$
generated by
all
$F_\alpha$,
$\alpha\in R^+$,
and
$K_i$,
$i\in[1, n]$,
each
$\lambda\in\Lambda$
defines a
1-dimensional
$\overline U_\zeta^\flat$-module
by annihilating
all
$F_\alpha$
and letting
$K_i$
act
by
$\zeta^{d_i\lambda_i}$.
Then
$\overline U_\zeta\otimes_{\overline U_\zeta^\flat}\lambda$
comes equipped with a structure of
$\tilde\fru$-module
[AJS, 2.10]
such that
each
$x\in U_\cK^0$
acts
on
$\prod_{\alpha\in R^+}
E_\alpha^{c_\alpha}\otimes1$,
$c_\alpha\in\bbN$,
by
the scalar
$\chi_{\lambda+\sum_\alpha
c_\alpha\alpha}(x)$.
Put
$\rho=\frac{1}{2}\sum_{\alpha\in R^+}\alpha$.

\begin{prop}
For each
$\lambda\in\Lambda$
we have an isomorphism of
$\tilde\fru$-modules
\[
\tilde\nabla(\lambda)\simeq
\overline U_\zeta\otimes_{\overline U_\zeta^\flat}
(\lambda-2(\ell-1)\rho).
\]
\end{prop}

\pf
Let
$\varepsilon_\lambda\in\tilde\nabla(\lambda)$
be the element
induced by the counit of
$\fru^+$.
By the universality of
$\tilde\nabla(\lambda)$
[APW, 0.8.1]
there is a homomorphism of
$\tilde\fru$-modules
\begin{equation}
\overline U_\zeta\otimes_{\overline U_\zeta^\flat}
(\lambda-2(\ell-1)\rho)\to
\tilde\nabla(\lambda)
\quad\text{such that}\quad
E^+\otimes1\mapsto
\varepsilon_\lambda,
\end{equation}
where
$E^+=\prod_{\alpha\in R^+}E_\alpha^{\ell-1}$.
On the other hand,
by
\cite[4.9]{AJS}
\begin{equation}
\text{
the
$\tilde\fru$-socle
of
$\overline U_\zeta\otimes_{\overline U_\zeta^\flat}(\lambda-2(\ell-1)\rho)$
is simple of highest weight
$\lambda$}.
\end{equation}
It follows that
the map (1) is injective,
and hence bijective by dimension.

\noindent
(1.5)
{\bf Corollary.}
{\it
Each
$\tilde\nabla(\lambda)$
(resp.
$\nabla(\lambda)$),
$\lambda\in\Lambda$,
is the projective cover of
$\lambda-2(\ell-1)\rho$
as
integrable
$\tilde\fru^\sharp$-
(resp.
$\fru^\sharp$-)
module.
}

\setcounter{equation}{0}
\noindent
(1.6)
Because of the isomorphism
(1.4)
we call
$\tilde\nabla(\lambda)$
and also by abuse of language
$\nabla(\lambda)$
the infinitesimal Verma module
of highest weight
$\lambda$.
By
\cite[6.3 and 4.10.1]{AJS}
\begin{equation}
\text{
$\tilde\nabla(\lambda)$
(resp.
$\nabla(\lambda)$)
is simple as
$\tilde\fru$-
(resp.
$\fru$-)
module
iff
$\lambda\equiv
(\ell-1)\rho\!\mod
\ell\Lambda$
}.
\end{equation}

Let
$\fru^{++}$
(resp. $\fru^{--}$)
be the augmentation ideal
of
$\fru^+$
(resp. $\fru^-$).
If
$M$
is a
$\fru^\pm$-module,
let
$M^{\fru^{\pm\pm}}$
denote the annihilater
of
$\fru^{\pm\pm}$
in
$M$.
By
(1.3)
\begin{equation}
\tilde\nabla(\lambda)^{\fru^{++}}
=
\nabla(\lambda)^{\fru^{++}}
=\lambda.
\end{equation}

If
$\lambda\equiv(\ell-1)\rho\mod
\ell\Lambda$,
then
$\nabla((\ell-1)\rho+\ell\nu)=
\nabla((\ell-1)\rho)$,
$\nu\in\Lambda$,
is simple,
called the Steinberg module,
and hence
\begin{equation}
\nabla((\ell-1)\rho+\ell\nu)^{\fru^{--}}
=
-(\ell-1)\rho+\ell\nu.
\end{equation}
In general,
the lowest weight of
$\tilde\nabla(\lambda)$
(resp. the socle of
$\tilde\nabla(\lambda)$)
is
$\lambda-2(\ell-1)\rho$
(resp.
$w_0\lambda^0+\ell\lambda^1$
if
$w_0$
is an element of the Weyl group of
$R$
such that
$w_0R^+=
-R^+$
and if
one writes
$\lambda=
\lambda^0+\ell\lambda^1$
with
$\lambda^0\in
\Lambda^{\res}$
and
$\lambda^1\in\Lambda$
[AJS, 4.2.5]).
It follows that
\begin{equation}
\text{
$\dim
\nabla(\lambda)^{\fru^{--}}\geq2$
unless
$\nabla(\lambda)$ is
$\fru$-simple}.
\end{equation}

\setcounter{equation}{0}
\noindent
(1.7)
Let
$\cB=\cA/(\phi_\ell)$
and let
$\tilde\fru_\cB$
be the $\cB$-subalgebra of
$U\otimes_\cA\cB$
generated by
$E_i\otimes1$,
$F_i\otimes1$,
$K_i\otimes1$,
$\lbr
K_i; c
\\
r
\rbr\otimes1$,
$i\in[1, n]$,
$c\in\bbZ$,
$r\in\bbN$.
Define
its $\cB$-subalgebras
$\tilde\fru_\cB$,
$\fru_\cB$,
$\tilde\fru_\cB^\flat$,
and
$\fru_\cB^\flat$
as for
$\tilde\fru$.
An infinitesimal Verma module may be defined over
$\cB$;
$\tilde\nabla_\cB(\lambda)=\tilde\fru_\cB^\flat\Mod(\tilde\fru_\cB,
\lambda)$,
$\lambda\in\Lambda$,
admits a structure of
$\tilde\fru_\cB$-module
like
$\tilde\nabla(\lambda)$,
and we have an isomorphism of
$\tilde\fru$-modules
\[
\tilde\nabla_\cB(\lambda)\otimes_\cB\cK\simeq
\tilde\nabla(\lambda).
\]
Restricting the
$\tilde\fru_\cB$-action to
$\fru_\cB$,
one obtains
$\fru_\cB$-module
$\nabla_\cB(\lambda)=
\fru_\cB^\flat\Mod(\fru_\cB,
\lambda)$.

Assume now
that
$\ell$
is a prime
$p$.
Then
$\cB/p\cB$
is a finite filed
$\bbF_p$
of $p$-elements.
Let
$G$ be a simply connected simple algebraic group over
$\bbF_p$
associated to the Cartan matrix
$A$
with a Borel subgroup
$B$ and a maximal torus
$T$ of
$B$
both split over
$\bbF_p$
such that
the roots of
$B$ are
$-R^+$.
Let
$G_1$
(resp.
$B_1$)
be the Frobenius kernel of
$G$
(resp.
$B$).
If
$\Dist(G_1)$
(resp.
$\Dist(B_1)$)
is the algebra of distibutions of
$G_1$
(resp.
$B_1$),
there are isomorphisms of
$\bbF_p$-algebras
[L2]
\[
\fru_\cB/(K_i-1\mid
i\in[1, n])\otimes_\cB\bbF_p
\simeq
\Dist(G_1),
\quad
\fru_\cB^\flat/(K_i-1\mid
i\in[1, n])\otimes_\cB\bbF_p
\simeq
\Dist(B_1),
\]
and
each
$\tilde\nabla_p(\lambda):=
\tilde\nabla_\cB(\lambda)\otimes_\cB\bbF_p$,
$\lambda\in\Lambda$,
admits a
structure of
$G_1T$-module
(cf. \cite[II.9]{J}):
\[
\tilde\nabla_p(\lambda)\simeq
\Dist(G_1)\otimes_{\Dist(B_1)}
(\lambda-2(p-1)\rho).
\]
Likewise
$\nabla_\cB(\lambda)\otimes_\cB\bbF_p$
yields a
$G_1$-module,
which we will denote by
$\nabla_p(\lambda)$.


\ \\
{\bf
$2^\circ$
Maximal cyclic modules}

In this section
we assume that our
Cartan matrix is of
type
$\rA_n$.
Then
all
$d_i=1$,
and we will suppress
$i$
from
$[\ \ ]_i$.
By properly
specializing the parametres
a maximal cyclic module
$\cV$
for
$U_\zeta$
of [DJMM]
factors through
$\fru$,
having a unique,
up to scalar,
$\fru^+$-primitive vector.
Thus
$\cV$
is of dimension
$\ell^{\vert
R^+\vert}$
and
has
a simple $\fru$-socle,
inviting us to
compare
$\cV$
with
infinitesimal
Verma modules.

\setcounter{equation}{0}
\ \\
(2.1)
Fix
$\lambda\in\Lambda$.
We define as we may
$\cV$
to be a
$\cK$-linear space of
basis
$u_\bm$,
$\bm\in
(\bbZ/\ell\bbZ)^{\vert
R^+\vert}$.
After \cite{N}
we
reindex
$R^+$
by the
pairs
$(i, j)$,
$1\leq
i\leq
j\leq
n$,
and we will denote the
$(i, j)$-component
of
$\bm$
by
$m_{ij}$.
Then
$\cV$
admits a structure
of
integrable
$\fru$-module
as proved in
\cite[5.2]{N}
such that
for each
$i\in[1, n]$
and
$\bm\in(\bbZ/\ell\bbZ)^{\vert R^+\vert}$
\begin{align}
E_iu_\bm
&=
\sum_{k=i}^n
[m_{ik}+
m_{i,k-1}-
m_{i-1,k-1}-
m_{i+1,k}]
u_{\bm+\epsilon(i,k)+\dots+\epsilon(i,n)},
\\
F_iu_\bm
&=
\sum_{k=1}^i
[-\lambda_i+
m_{i+1-k,n-k}-
m_{i+1-k,n+1-k}+
m_{i-k,n+1-k}-
m_{i-k,n-k}]
\\
&\hspace*{1cm}
u_{\bm-\epsilon(i+1-k,n+1-k)-\epsilon(i+2-k,n+2-k)-\dots-\epsilon(i,n)},
\notag
\\
K_iu_\bm
&=
\zeta^{
\lambda_i+
2m_{in}-
m_{i-1,n}-
m_{i+1,n}}u_\bm,
\end{align}
where
$[r]=
\frac{\zeta^r-\zeta^{-r}}{\zeta-\zeta^{-1}}$,
$\epsilon(i, j)\in(\bbZ/\ell\bbZ)^{\vert R^+\vert}$
such that
$\epsilon(i, j)_{st}=
\delta_{is}\delta_{jt}$
for each
$s$
and
$t$,
and
any meaningless
terms in the sums
should be read as $0$.
As the structure of
$\fru$-module on
$\cV$ depends on
$\lambda$,
to be precise, we will denote the
$\fru$-module $\cV$
by
$\cV(\lambda)$.

A main theorem of [N]
is that
$\cV(\lambda)$ has a unique,
up to $\cK^\times$,
$\fru^+$-primitive element,
i.e.,
\begin{equation}
\cV(\lambda)^{\fru^{++}}=
\cK
u_{\vec{0}}
\quad\text{with}
\quad
\vec{0}=(0,
\dots, 0),
\end{equation}
and hence by Engel's theorem
\begin{equation}
\text{
$\cV(\lambda)$
has a simple
$\fru$-socle
generated by
$u_{\vec{0}}$}.
\end{equation}
It
also
follows from (1.3)
by dimension that
\begin{prop}
There is an isomorphism of
$\fru^\sharp$-modules
$\cV(\lambda)\simeq\nabla(\lambda)$.
\end{prop}

\setcounter{equation}{0}
\noindent
(2.2)
Recall anitiautomorphism
$\tau$ of
$\fru$ such that
\[
E_i\mapsto
F_i,
\quad
F_i\mapsto
E_i,
\quad
K_i\mapsto
K_i,
\quad\forall
i\in[1, n].
\]
If
$M$ is
a
$\fru$-module,
let
$M^\tau$
be the
$\cK$-dual space of
$M$ with
a
$\fru$-action given by
\[
xf=f(\tau(x)?),
\quad
f\in M^*,
x\in\fru.
\]
Then the isomorphism of
$\fru^\sharp$-modules
$\cV(\lambda)\simeq
\nabla(\lambda)$
from
(2.1) yields
an isomorphism
of
$\fru^-$-modules
\begin{align*}
\cV(\lambda)^\tau
&\simeq
\nabla(\lambda)^\tau
\\
&\simeq
\{
\overline{U}_\zeta\otimes_{\overline{U}_\zeta^\flat}
(\lambda-2(\ell-1)\rho)\}^\tau
\quad\text{by (1.4)}
\\
&\simeq
\overline{U}_\zeta\otimes_{\overline{U}_\zeta^\sharp}\lambda
\quad\text{by [AJS, 4.10]}
\\
&\simeq
\fru^-,
\end{align*}
where
$\overline{U}_\zeta^\sharp$
is the
$\cK$-subalgebra of
$\overline{U}_\zeta$
generated by
all
$E_\alpha$,
$\alpha\in R^+$,
and
$K_i$,
$i\in[1, n]$.
It follows from
[X, 2.5] again
that
$\cV(\lambda)^\tau$
has a unique,
up to
$\cK^\times$,
$\fru^-$-primitive vector,
and hence
\begin{cor}
$\cV(\lambda)$
has the same simple $\fru$-socle and
the same simple $\fru$-head head
as
$\nabla(\lambda)$.
\end{cor}

\setcounter{equation}{0}
\noindent
(2.3)
We find, moreover, that
\begin{thm}
The
$\fru$-module
$\cV(\lambda)$
has a unique,
up to
$\cK^\times$,
$\fru^-$-primitive vector
$u_\bm$ with
\[
m_{ij}\equiv
-\sum_{s=1}^i
\sum_{t=i}^j
\lambda_{n+s-t}
\mod
\ell
\quad\forall
i\leq
j,
\]
which has
$\fru^0$-weight
$w_0\lambda$.
\end{thm}

\pf
The argument is the same as for
(2.1.4) from [N, 4.2];
let
\[
\sum_{\bm\in(\bbZ/\ell\bbZ)^{\vert
R^+\vert}}
c_\bm
u_\bm\in
\cV(\lambda)^{\fru^{--}},
\quad
c_\bm\in\cK.
\]
As
$F_i\sum
c_\bm
u_\bm=0$
for all
$F_i$,
we obtain successively
$c_\bm=0$
unless
\begin{equation}
-\lambda_i+
m_{i+1-k,n-k}-m_{i+1-k,n+1-k}+
m_{i-k,n+1-k}-m_{i-k,n-k}
\equiv0\mod\ell
\end{equation}
for each
$i$
and
$k\in[1, i]$.
If
$\sum
c_\bm
u_\bm\ne0$,
the system (1)
of equations
determines
$\bm$
with
$c_\bm\ne0$
uniquely
as asserted.

\noindent
(2.4)
{\bf Corollary.}
{\it
Let
$\lambda\in X$.
\begin{enumerate}
\item[{\rm (i)}]
There is an isomorphism of
$\fru^\flat$-modules
\[
\cV(\lambda)\simeq
\overline{U}_\zeta\otimes_{\overline{U}_\zeta^\sharp}
(w_0\lambda+2(\ell-1)\rho).
\]
\item[{\rm (ii)}]
The
$\fru$-module
$\cV(\lambda)$
lifts to
an integrable
$\tilde\fru$-module
iff
$\lambda\equiv(\ell-1)\rho\mod
\ell\Lambda$.
In particular,
unless
$\lambda\equiv(\ell-1)\rho\mod
\ell\Lambda$,
$\cV(\lambda)$
is not
isomorphic as
$\fru$-module
to
any infinitesimal Verma module.
\end{enumerate}
}

\pf
For
(i) argue as in (2.1) and (2.2).

(ii)
If
$\lambda\equiv(\ell-1)\rho\mod
\ell\Lambda$,
the simple
$\fru$-socle
of
$\cV(\lambda)$
has
dimension
$\ell^{\vert
R^+\vert}=\dim
\cV(\lambda)$
and hence the assertion follows.
We may therefore assume that
$\lambda\not\equiv(\ell-1)\rho\mod
\ell\Lambda$.

Just suppose
$\cV(\lambda)$ lift to
an integrable
$\tilde\fru$-module.
Then we would have from (1.3)
an isomorphism
of
$\tilde\fru^\sharp$-modules
\[
\cV(\lambda)\simeq
\tilde\nabla(\lambda+\ell\nu)
\quad\text{
for some
$\nu\in\Lambda$}.
\]
Then the unique $\fru^-$-primitive
in $\cV(\lambda)$
should have by (1.4)
weight
$\lambda+\ell\nu-2(\ell-1)\rho$.
That
would yield,
arguing as in (2.3),
an isomorphism of
$\tilde\fru^\flat$-modules
\[
\cV(\lambda)\simeq
\overline U_\zeta\otimes_{\overline U_\zeta^\sharp}
(\lambda+\ell\nu).
\]
Then
$\cV(\lambda)=\fru^-u_{\vec{0}}$
as
$u_{\vec{0}}$
is a highest weight vector of
$\cV(\lambda)$.
But
$\fru^-u_{\vec{0}}=
\fru
u_{\vec{0}}$
is
the simple socle of
$\cV(\lambda)$
and of dimension
$<
\dim\cV(\lambda)$
by
(1.6.1),
absurd.

\setcounter{equation}{0}
\noindent
(2.5)
Assume now
that
$\ell$
is a prime
$p$
and
let
$\cB=\cA/(\phi_p)=
\bbZ[v]/(\phi_p)$
as in (1.7).
By construction
$\cV(\lambda)$
may be defined over
$\cB$:
let
$\cV_\cB(\lambda)$
be the free
$\cB$-module
of basis
$u_\bm$,
$\bm\in(\bbZ/p\bbZ)^{\vert
R^+\vert}$,
with
the $\fru_\cB$-action
given by
(2.1.1-3).
Regarding
$\bbF_p$
as the quotient
$\cB/p\cB=\cB/(v-1)$,
let
$\cV_p(\lambda)=\cV_\cB(\lambda)\otimes_\cB\bbF_p$.
Then
$\cV_p(\lambda)$
is naturally a
$G_1$-module
for
$G=\SL_{n+1}\otimes_\bbZ\bbF_p$
in the setup of
(1.7).
The proof of
(2.1.4)
from
\cite{N}
and the
argument of (2.3)
carry over
to obtain
\begin{thm}
Let
$\lambda\in\Lambda$
\begin{enumerate}
\item[\rm{(i)}]
There is an isomorphism of
$B^+_1$-modules
$
\cV_p(\lambda)\simeq
\nabla_p(\lambda)$.
\item[\rm{(ii)}]
$\cV_p(\lambda)$
has a
unique,
up to
$\bbF_p^\times$,
$B^+_1$-
and
$B_1$-invariant vector,
respectively,
and hence has
a
simple
$G_1$-socle
and a
simple
$G_1$-head,
where
$B_1^+$
is the Frobenius kernel
of the
Borel subgroup
$B^+$
of $G$
opposite to
$B$.
\item[\rm{(iii)}]
If
$\lambda\not
\equiv
(p-1)\rho
\mod
p\Lambda$,
the structure of
$G_1$-module
on
$\cV_p(\lambda)$
does not
lift to
$G_1T$-module,
and
hence
$\cV_p(\lambda)$
is not isomorphic to
any infinitesimal
Verma module
as
$G_1$-module.
\end{enumerate}
\end{thm}

\noindent
(2.6)
{\bf Remark.}
If
the simple
$\fru$-module
generated by
$u_{\vec{0}}$
in
$\cV(\lambda)$
and
the simple
$G_1$-module
generated by
$u_{\vec{0}}\otimes1$
in
$\cV_\cB(\lambda)\otimes_\cB\bbF_p$
have the same dimension,
then
Lusztig's conjecture for the irreducible
$\SL_{n+1}\otimes_\bbZ\bbF_p$-modules
will follow from
\cite{KL}
and
\cite{KT},
\cite{C}.



\begin{thebibliography}{AAAA}

\bibitem[AJS]{AJS}
Andersen, H. H.,
Jantzen, J. C.
and Soergel, W.,
{\it Representations of quantum groups at a p-th root of unity and of
semisimple
groups in characteristics p: Independence of p}, Ast{\'e}risque {\bf 220} (
1994), 1-320.
\bibitem[APW]{APW}
Andersen, H. H.,
Polo, P.
and Wen, K., {\it Injective modules for
quantum algebras}, Am.  J. Math. {\bf 114} (1992), 571-604.
\bibitem[C]{C} Casian, L.,
{\it Proof of the
Kazhdan-Lusztig conjecture for Kac-Moody algebras
(the characters
$\ch
L_{w\rho-\rho}$)},
Adv. Math. {\bf 119} (1996), 207--281.
\bibitem[DJMM]{DJMM} Date, E.,
Jimbo, M.,
Miki, K. and Miwa, T., {\it Cyclic representations of
$U_q(\mathfrak{sl}(n+1, \bbC))$
at
$q^N=1$},
Publ. RIMS
{\bf 27} (1991), 347--366.
\bibitem[DK]{DK} De Concini, C.
and
Kac, V.,
{\it Representations of quantum groups at roots of 1}, pp.
471--506
in:
A. Conne et al.
(eds.),
{\it
Operator Algebras,
Unitary Representations,
Enveloping Algebras, and Invariant Theory}
(Colloque Dixmier),
Proc. Paris 1989
(Progress in Mathematics
{\bf 92}),
Boston etc. 1990
(Birkh\"{a}user)
\bibitem[H]{H}
Humphreys, J. E.,
{\it Modular representations
of classical Lie algebras and semisimple groups}, J.
Alg. {\bf 19} (1971), 51--79.
\bibitem[J]{J} Jantzen, J. C.,
\newblock \textsl{Representations of
Algebraic Groups},
\newblock Pure and Applied Mathematics {\bf 131},
\newblock Academic Press, Boston etc. 1987.
\bibitem[KT]{KT}
Kashiwara, M. and Tanisaki, T.,
{\it Kazhdan-Lusztig conjecture
for affine
Lie algebras with negative level},
Duke Math. J.
{\bf 77}
(1995),
21--62
\bibitem[KL]{KL}
Kazhdan, D.A. and
Lusztig, G.,
{\it Tensor structures arising
from affine Lie algebras
I, II},
J. AMS
{\bf 6} (1993),
905--1011;
{\it Tensor structures arising from affine Lie algebras
III, IV},
J. AMS
{\bf 7} (1994),
335--453.
\bibitem[L1]{L1} Lusztig, L.,
{\it Finite dimensional Hopf algebras arising from quantized universal
enveloping algebras}, J. AMS. {\bf 3} (1990), 257--296.
\bibitem[L2]{L2} Lusztig, L.,
{\it Quantum groups at roots of 1}, Geom. Ded. {\bf 35} (1990), 89--114.
\bibitem[N]{N} Nakashima, T.,
{\it Irreducible modules of finite dimensional quantum algebras of type $A$
at roots of unity},
J. Math. Phys. {\bf 43} (2002), no. 4, 2000-2014.
\bibitem[X]{X} Xi, N.,
{\it Irreducible modules of
quantized enveloping algebras at roots of
1},
Publ. RIMS
{\bf 32} (1996), 235--276.












\end{thebibliography}
\end{document}